\documentstyle[12pt,amssymb]{amsart}
\input{diagrams.tex}
\newcommand{\PP}{{\mathbb P}}
\newcommand{\T}{{\mathbb T}}

\newcommand{\C}{{\mathbb C}}

\renewcommand{\phi}{\varphi}

\title{\bf Hyperbolicity of general deformations}

\author{\bf Mikhail Zaidenberg}

\address{Universit\'e
Grenoble I, Institut Fourier, UMR 5582 CNRS-UJF, BP 74, 38402 St.\
Martin d'H\`eres c\'edex, France}
\email{zaidenbe@@ujf-grenoble.fr}

\thanks{
\mbox{\hspace{11pt}}{\it 2000 Mathematics Subject Classification}:
14J70,  32J25.\\
\mbox{\hspace{11pt}}{\it Key words}: Kobayashi hyperbolicity,
projective hypersurface, deformation}

\date{}

\begin{document}
\begin{abstract} This is the content of the talk
given at the conference ``Effective Aspects of Complex Hyperbolic
Varieties", Aver Wrac'h, France, September 10-14, $'$07. We
present two methods of constructing low degree Kobayashi
hyperbolic hypersurfaces in $\PP^n$:

\smallskip

\noindent {$\bullet$ The projection method}

\smallskip

\noindent {$\bullet$ The deformation method}

\smallskip

\noindent The talk is based on joint works of the speaker with B.\
Shiffman and C.\ Ciliberto.
\end{abstract}

\maketitle

\vskip 0.2in

\section{\bf DIGEST on KOBAYASHI THEORY}

\subsection{Kobayashi hyperbolicity}$\,$

\vskip 0.2in

{  \noindent {\bf DEFINITION} The {\it Kobayashi pseudometric}
$k_X$ on a complex space $X$ satisfies the following axioms :

\medskip

\noindent (i) On the unit disc $\Delta$, the Kobayashi
pseudometric $k_\Delta$ coincides with the Poincar\'e metric;

\medskip

\noindent (ii) every holomorphic map $\varphi:\Delta \to X$ is a
contraction: $\varphi^*(k_X)\le k_\Delta$;

\medskip

\noindent (iii) $k_X$ is the maximal pseudometric on $X$
satisfying (i) and (ii).

\vskip 0.2in

\noindent {\bf REMARK} Every holomorphic map $\varphi: X \to Y$ is
a contraction: $\varphi^*(k_Y)\le k_X$.

\vskip 0.2in

\noindent {\bf DEFINITION} $X$ is called {\it Kobayashi
hyperbolic} if $k_X$ is non-degenerate : $$k_X(p,\,q)=0
\Longleftrightarrow p=q.$$

\vskip 0.2in

\noindent {\bf EXAMPLES} $k_{\C^n}\equiv 0,\quad k_{\PP^n}\equiv
0,\quad k_{\T^n}\equiv 0\,,$ where $\T^n=\C^n/\Lambda$ is a
complex torus,

\vskip 0.2in

\noindent whereas $\C\setminus \{0,\,1\}$ is hyperbolic (the
Schottky-Landau Theorem.)

\vskip 0.2in


\subsection{\bf Classical theorems}$\,$

\vskip 0.2in

\noindent According to the above definition and to Royden's
Theorem, $X$ is hyperbolic iff natural analogs of the classical
Schottky and Landau Theorems hold for $X$.

\vskip 0.2in

\noindent {\bf Brody-Kiernan-Kobayashi-Kwack THEOREM}
\\ {\it For a compact complex space $X$ the following conditions
are equivalent :

\medskip

\noindent $\bullet$ $X$ is Kobayashi hyperbolic;

\medskip

\noindent $\bullet$ Little Picard Theorem holds for $X$ :
$$\forall f\,:\,\C \to X,\,\,\,\,f= {\rm const};$$
\noindent $\bullet$  Big Picard Theorem holds for $X$ :
$$\forall f\,:\,\Delta\setminus\{0\} \to X\,\,\,\,\exists
{\bar f}\,:\,\Delta \to X,\,:\,{\bar f}\vert
(\Delta\setminus\{0\})=f;$$ \noindent $\bullet$ Montel
 Theorem holds for $X$ :
the space $HOL(\Delta,\,X)$ is compact.}

\vskip 0.2in

\noindent {\bf REMARK } If $X$ is hyperbolic then $\forall Y$, the
space  $HOL(Y,\,X)$ is compact.

\vskip 0.2in

\noindent {\bf DEFINITION} Let $M$ be  a hermitian compact complex
manifold.  An entire curve $\varphi : \C\to M$ is called a {\it
Brody curve} if $$||\varphi'(z)||\le 1=||\varphi'(0)||\quad
\forall z\in\C\,.$$

\vskip 0.2in

\noindent {\bf Brody's THEOREM} {\it $M$ as above is hyperbolic
iff it does not possess any Brody entire curve.}

\vskip 0.2in

\noindent {\bf Brody's STABILITY THEOREM} \\{\it Every compact
hyperbolic subspace $X$ of a complex space $Z$ admits a hyperbolic
neighborhood. Consequently, every compact subspace $X'\subseteq Z$
sufficiently close to $X$ is hyperbolic. In particular, if
$X\subseteq \PP^n$ is a hyperbolic hypersurface then every
hypersurface $X'\subseteq \PP^n$ sufficiently close to $X$ is
hyperbolic too.}

\vskip 0.2in

\subsection{Hyperbolicity of hypersurfaces in $\PP^n$}$\,$

\vskip 0.2in

\centerline{\underline{\bf Kobayashi Problem ($'$70)}}

\vskip 0.2in

\noindent {\it Is it true that a (very) general hypersurface $X$
of degree $d\ge 2n-1$ in $\PP^n$ is Kobayashi hyperbolic? In
particular, is this true for a (very) general surface $X$ in
$\PP^3$ of degree $d\ge 5$ ?}

\vskip 0.2in

 \centerline{\underline{\bf Hyperbolic surfaces in $\PP^3$}}

\vskip 0.2in

\noindent {\bf THEOREM (McQuillen \cite{MQ}, Demailly-El Goul \cite{DEG})}\\
{\it A very general surface $X$ in $\PP^3$ of degree $d\ge 21$ is
Kobayashi hyperbolic.}

\vskip 0.2in

\noindent For some recent advances in higher dimensions, see
Y.-T.\ Siu \cite{Si} and E.\ Rousseau \cite{Ro}.

\vskip 0.2in

\centerline{\underline{\bf EXAMPLES}} \centerline{\underline{\bf
of small degree hyperbolic surfaces in $\PP^3$}}

\vskip 0.2in

\noindent Concrete examples were found by \\
{\bf Brody-Green} $'$77, $d=2k\ge 50$,\\ {\bf Masuda-Noguchi}
$'$96, $d=3e\ge 24$,\\ {\bf Khoai} $'$96, $d\ge 22$,\\ {\bf Nadel}
$'$89, $d\ge 21$,\\ {\bf Shiffman-Z$'$} $'$00, $d\ge 16$,
\\{\bf El Goul} $'$96, $d\ge 14$,
\\{\bf Siu-Yeung} $'$96, {\bf Demailly-El Goul} $'$97, $d\ge 11$,\\
{\bf J.\ Duval} $'$99 \cite{Du2}, {\bf Shirosaki-Fujimoto} $'$00
\cite{Fu}, $d=2k\ge 8$:} \begin{equation}\label{1}
Q(X_0,X_1,X_2)^2-P(X_2,X_3)=0\,,\end{equation} where $Q,\,P$ are
generic homogeneous formes of degree $k$ and
$d=2k$, respectively,\\
{\bf Shiffman-Z$'$} $'$02 \cite{SZ}, $d=8$,\\
{\bf Shiffman-Z$'$}
 $'$05 \cite{SZ1}, $d\ge 8$,\\
{\bf J.\ Duval} $'$04 \cite{Du}, $d=6$.

\vskip 0.2in

\noindent Algebraic families of hyperbolic hypersurfaces
$X_n\subseteq\PP^n$ for any
$n\ge 3$ were constructed e.g., by \\
{\bf Masuda-Noguchi} $'$96,\\
{\bf Siu-Yeung} $'$97,\\
{\bf Shiffman-Z$'$} $'$02 \cite{SZ2}.

\vskip 0.2in

\noindent In these examples $\deg X_n$ grows quadratically with
$n$, for instance, $\deg X_n=4(n-1)^2$ \cite{SZ2}. Whereas the
Kobayashi Conjecture suggests a linear growth of the minimal such
degree. This leads to the following problem.

\vskip 0.2in

\noindent {\bf PROBLEM} {\it Find a sequence of hyperbolic
hypersurfaces $X_n\subseteq\PP^n$ with $\deg X_n\le Cn$ for some
positive constant $C$.}

\section{\bf PROJECTION METHOD}

\bigskip

\subsection{\bf Symmetric powers of curves as hyperbolic
hypersurfaces}$\,$

\vskip 0.2in

\noindent {\bf PROPOSITION (Shiffman-Z$'$ $'$00 \cite{SZ3})} {\it
The $n$th symmetric power $C^{(n)}$ of a generic smooth projective
curve $C$ of genus $g\ge 3$ is hyperbolic iff $g\ge 2n-1$. In
particular, the symmetric square $C^{(2)}$ is always hyperbolic.}

\vskip 0.2in

\noindent {\bf THEOREM (Shiffman-Z$'$ $'$00 \cite{SZ3})} {\it With
$C$ as before, let us consider an embedding
$C^{(2)}\hookrightarrow\PP^5$. Then a general projection  $S$ of
$C^{(2)}$ to $\PP^3$ is hyperbolic. The minimal degree of such a
hyperbolic surface $S\subseteq\PP^3$ is equal $16$.}

\vskip 0.2in

\noindent {\bf EXAMPLE of degree 16:} Let $C\subseteq\PP^2$ :
$x^4-xz^3-y^3z=0\,,$ and let $C^{(2)}\hookrightarrow\PP^5$ be
embedded via the natural embedding of the symmetric square of
$\PP^2$ in $\PP^5$. Then a general projection  of $C^{(2)}$ to
$\PP^3$ is a singular hyperbolic surface $S\subseteq \PP^3$ of
degree $16$, with the double curve $D$ of genus $142$.

\vskip 0.2in

\noindent Let us explain in brief our methods. Let
$V\hookrightarrow\PP^5$ be a smooth hyperbolic surface, and let
$\pi:V\to S\hookrightarrow\PP^3$ be a projection. Then $S$ has
self-intersection along a double curve $D\subseteq S$. By the
universal property of the normalization, there is a commutative
diagram

\begin{diagram}
V&&\rTo^{\psi}&&S_{\rm norm}\\
&\rdTo_\pi & &\ldTo_\nu& \\
&&S&&
\end{diagram}

\noindent where $\nu:S_{\rm norm}\to S$ is the normalization. By
Zariski's Main Theorem, $\psi:V\to S_{\rm norm}$ is an
isomorphism. Hence any entire curve $\varphi:\C\to S$ can be
lifted to $V=S_{\rm norm}$:
\begin{diagram}
&&V\\
&\ruTo^{\tilde{\varphi}} &\dTo_\nu& \\
\C&\rTo_\varphi&S&&
\end{diagram}
unless $\varphi(\C)\subseteq D$. Since $V$ is hyperbolic,
${\tilde{\varphi}}=$cst. Thus $S$ is hyperbolic iff $D$ is. A
similar argument shows that $S$ is always hyperbolic modulo $D$.
In the proof of the above theorem we show that, for a general
projection, $D$ is hyperbolic indeed and so $S$ is. Similarly, for
the Cartesian square of a curve the following holds.

\vskip 0.2in

\noindent {\bf PROPOSITION (Shiffman-Z$'$ $'$00 \cite{SZ3})} {\it
Let $C$ be a smooth projective curve of genus $g\ge 2$. Let us fix
an embedding $V=C\times C\hookrightarrow\PP^n\,.$ Then the double
curve $D\subseteq S$ of a general projection $V\to S\subseteq
\PP^3$ is irreducible of genus $g(D)\ge 225$, and $S$ is a
singular hyperbolic surface of degree $\ge 32$.}

\vskip 0.2in

\noindent However for a non-generic projection, the double curve
of the image surface can be neither irreducible nor hyperbolic.

\vskip 0.2in

\noindent {\bf EXAMPLE (Kaliman-Z$'$ $'$01 \cite{KZ})} Consider
the smooth Fermat quartic $$C:\,x^4+y^4+z^4=0\quad\mbox{ in}\quad
\PP^2\,.$$ Then the product $V=C\times C$ admits a projective
embedding and a projection to $\PP^3$ such that the double curve
$D$ of the image surface $S\subseteq \PP^3$ consists of  $4$
disjoint projective lines. Thus $S$ is not hyperbolic whereas its
normalization $V$ is.

\vskip 0.2in

\noindent For 3-folds in $\PP^4$ we have the following result.

\vskip 0.2in

\noindent {\bf THEOREM (Ciliberto-Z$'$ $'$03 \cite{CZ})} {\it For
a general projective curve $C$ of genus $g\ge 7$, we fix an
embedding
 $C^{(3)}\hookrightarrow\PP^7$. Then a general projection $X$ of
$C^{(3)}$ to $\PP^4$ is a hyperbolic hypersurface in $\PP^4$. This
is also true for a general quintic $C\subseteq \PP^2$ ($g=6$) and
a certain special embedding $C^{(3)}\hookrightarrow\PP^7$ of
degree $125$. The latter is the minimal degree which can be
achieved via the projection method using the symmetric cubes
$C^{(3)}$.}

\vskip 0.2in

\noindent The proof goes as follows. It is shown that

\vskip 0.2in

\noindent $\bullet$ $C^{(3)}$ does not contain any curve of genus
$<g$; in particular, it is hyperbolic.

\vskip 0.2in

\noindent $\bullet$ $X\subseteq\PP^4$ is hyperbolic iff the double
surface $S={\rm sing}\,(X)$ is. This uses the above trick with
lifting entire curves to the normalization $C^{(3)}$ of $X$.

\vskip 0.2in

\noindent $\bullet$ The irregularity $q(S)\ge g>5$. This is based
on  the fact that for a curve $C$ with general moduli, the
Jacobian $J(C)$ is a simple abelian variety.

\vskip 0.2in

\noindent $\bullet$ $S$ is hyperbolic iff it is algebraically
hyperbolic that is, does not contain any rational or elliptic
curve. This is based on the Bloch Conjecture.

\vskip 0.2in

\noindent $\bullet$ $S$ is hyperbolic iff the triple curve
$T\subseteq S$ of $X$ is. Recall that in a general point of $T$, 3
smooth branches of $X$ meet transversally. Actually $T$
parameterizes the 3-secant lines of $C^{(3)}\subseteq\PP^7$
parallel to the center of the projection
$\PP^7\dashrightarrow\PP^4$. The proof is based on Pirola's and
Ciliberto-van der Geer's results on deformations of hyperelliptic
and bielliptic curves on abelian varieties.

\vskip 0.2in

\noindent $\bullet$ Any irreducible component of the triple
curve\footnote{Presumably $T$ is irreducible, but we don't dispose
a proof of this.} $T$ has genus $\ge 2$. The proof is rather
involved.

\vskip 0.5in

\section{\bf DEFORMATION METHOD}

\vskip 0.2in

\noindent Let $X_0=f_0^*(0),\, X_{\infty}=f_\infty^*(0)$ be two
hypersurfaces of the same degree $d$ in $\PP^n$, and let
$$\{X_t\}_{t\in \PP^1}=\langle X_0, X_{\infty}\rangle,\qquad\mbox
{where}\qquad X_t=(f_0+t f_\infty)^*(0)\,,$$ be the pencil of
hypersurfaces generated by $X_0$ and $X_{\infty}$. For small
enough $|\varepsilon|\neq 0$ we call $X_\varepsilon$ a small
(linear) deformation of $X_0$ in direction of $X_\infty$.

\vskip 0.2in

\noindent {\bf DEFINITION } We say that a (very) general small
deformation of $X_0$ is hyperbolic if $X_\varepsilon$ is for a
(very) general $X_{\infty}$ and for all sufficiently small
$\varepsilon\neq 0$ (depending on $X_{\infty}$).

\vskip 0.2in

Let us formulate the following

\vskip 0.2in

\noindent  {\bf ``Weak Kobayashi Conjecture" :} {\it For every
hypersurface $X\subseteq\PP^n$ of degree $d\ge 2n-1$, a (very)
general small deformation of $X$ is Kobayashi hyperbolic.}

\vskip 0.2in

\noindent By Brody's Theorem, the proof of hyperbolicity of $X$
reduces to a certain degeneration principle for entire curves in
$X$. The Green-Griffiths' 79$'$ proof of Bloch's Conjecture
\cite{GG} provides a kind of such degeneration principle. It was
shown by McQuillen \cite{MQ} and, independently, by Demailly-El
Goul \cite{DEG} (according with this principle) that every entire
curve $\varphi:\C\to X$ in a very general surface
$X\subseteq\PP^3$ of degree $d\ge 36$ ($d\ge 21$, respectively)
satisfies a certain algebraic differential equation.

\vskip 0.2in

\noindent Consider again a pencil $(X_t)$. Assuming that for a
sequence
 $\varepsilon_n\to 0$ the hypersurfaces $X_{\varepsilon_n}$
 are not hyperbolic, one can find a sequence
 of Brody entire curves
 $\varphi_n:\C\to X_{\varepsilon_n}$ which
 converges to a (non-constant) Brody curve
 $\varphi:\C\to X_{0}$.

\vskip 0.2in

\noindent Suppose in addition that
 $X_{0}$ admits
 a rational map $\pi:X_{0}\dashrightarrow Y_0$
 to a hyperbolic variety $Y_0$
 (to a curve $Y_0$ of genus $\ge 2$ in case
$\dim X_{0}=2$).
 Then necessarily
 $\pi\circ \varphi=\,$cst, provided that
 the composition $\pi\circ \varphi$
 is well defined. Anyhow
 the limiting Brody curve
 $\varphi:\C\to X_{0}$ degenerates.
 This degeneration however is
 not related to any specific property
 of the configuration $X_0\cup X_\infty$, but of
 $X_0$ alone. Here is another degeneration principle which
 involves both $X_0$ and $X_\infty$.

\vskip 0.2in

\noindent {\bf PROPOSITION 1 (Shiffman-Z$'$ $'$05 \cite{SZ}, Z$'$
$'$07 \cite{Za})} Consider a pencil of degree $d$ hypersurfaces
$X_\varepsilon \subseteq\PP^{n+1}$ generated by $X_0=X_0'\cup
X_0''$ and $ X_\infty$. Let $D=X_0'\cap X_0''$. Then for any
sequence of entire curves $\varphi_{n}:\C\to X_{\varepsilon_n}$
which converges to $\varphi:\C\to X_0'$ the following alternative
holds:\\ $\bullet$ Either $\varphi (\C)\subseteq D$, or\\
$\bullet$ $\varphi (\C)\cap D\subseteq D\cap X_\infty$ and
$d\varphi (t)\in T_P X_0'\cap T_P X_\infty\quad\forall P=\varphi
(t)\in D\cap X_\infty$.

\vskip 0.2in

\noindent {\bf THEOREM 1 (Z$'$ $'$07 \cite{Za})} {\it Let $Y_0$ be
a Kobayashi hyperbolic hypersurface  of degree $d$ in $\PP^n$
($n\ge 2$), where $\PP^n$ is realized as the hyperplane
$H=\{z_{n+1}=0\}$ in $\PP^{n+1}$. Then a general small deformation
 $X_\varepsilon\subseteq\PP^{n+1}$ of the double cone $2X_0$ over $Y_0$
is Kobayashi hyperbolic.}

\vskip 0.2in

\noindent The proof is based on Proposition 1 and on the following
lemma.

\vskip 0.2in

\noindent {\bf LEMMA 1} {\it Let $\hat Y\subseteq\PP^{n+1}$ be a
cone over a projective variety $Y\subseteq\PP^n$, and let
$X'\subseteq \PP^{n+1}$ be a general hypersurface  of degree $e\ge
2\dim Y$. Then $X'$ meets every generator $l$ of $\hat Y$ in at
least $k=e-2\dim Y$ points transversally.}

\vskip 0.2in

\noindent {\bf Proof of Theorem 1.} Suppose the contrary. Then we
can find a sequence $\varepsilon_n\longrightarrow 0$ and a
sequence of Brody curves $\varphi_n:\C\to X_{\varepsilon_n}$ such
that $\varphi_n\longrightarrow \varphi$, where $\varphi:\C\to X_0$
is non-constant. We let $\pi:X_0\dashrightarrow Y_0$ be the cone
projection. Since $Y_0$ is assumed to be hyperbolic we have
$\pi\circ\varphi=\,\,$cst. In other words $\varphi (\C)\subseteq
l$, where $l\cong\PP^1$ is a generator of the cone $X_0$.

\smallskip

\noindent We note that $\bigtriangledown f_0^2 |_{X_0}=0$. If $l$
and $X_\infty$ meet transversally in a point $\varphi (t)\in l\cap
X_\infty$ then $d\varphi(t)=0$ by virtue of
Proposition 1.\\


\noindent Since $Y_0\subseteq \PP^n$ is hyperbolic and $n\ge 2$ we
have $d\ge n+2$. In particular $$\deg X_\infty=2d\ge 2n+4\ge 2\dim
Y +5\,.$$ By Lemma 1, $l$ and $X_\infty$ meet transversally in at
least 5 points. Hence the nonconstant meromorphic function
$\varphi:\C\to l\cong\PP^1$ possesses at least 5 multiple values.
Since the defect of a multiple value is $\ge 1/2$, this
contradicts the Defect Relation. \qed

\vskip 0.2in

\noindent {\bf REMARK} Given a hyperbolic hypersurface $Y\subseteq
\PP^n$ of degree $d$, Theorem 1 provides a hyperbolic hypersurface
$X\subseteq \PP^{n+1}$ of degree $2d$. Iterating the construction
yields hyperbolic hypersurfaces in $\PP^n$ $\forall n\ge 3$ of
degree that grows exponentially with $n$.

\vskip 0.2in

\noindent {\bf EXAMPLE (Z$'$ $'$07 \cite{Za})} Let $C\subseteq
\PP^2$ be a hyperbolic curve of degree $d\ge 4$, and let
$X_0\subseteq \PP^3$ be a cone over $C$. Then a general small
deformation of the double cone $2X_0$ is a Kobayashi hyperbolic
surface in $\PP^3$ of even degree $2d\ge 8$.

\vskip 0.2in

\noindent The following example combines the projection and
deformation methods.

\vskip 0.2in

\noindent {\bf EXAMPLE (Shiffman-Z$'$ $'$03 \cite{SZ1})} There is
a singular octic $X_0\subseteq\PP^3$ whose normalization is a
simple abelian surface. Moreover, a general small deformation of
$X_0$ is Kobayashi hyperbolic.

\vskip 0.2in

\noindent {\bf EXAMPLE (Shiffman-Z$'$ $'$05 \cite{SZ})} Let
$X_0=X_0'\cup X_0''$ be the union of two cones in general position
in $\PP^3$ over smooth plane quartics $C',\,C''\subseteq\PP^2$,
respectively. Then a general small deformation of $X_0$ is
Kobayashi hyperbolic.

\vskip 0.2in

\noindent {\bf Sketch of the proof.} Suppose that for a sequence
$\varepsilon_n\to 0$, $X_{\varepsilon_n}$ is not hyperbolic. Then
we can find a sequence of Brody curves $\varphi_{n}:\C\to
X_{\varepsilon_n}$ which converges to a Brody curve $\varphi:\C\to
X_0$. We may assume that $\varphi(\C)\subseteq X_0'$.

Since $C'$ has genus $3$, $\pi'\circ \varphi : \C \to C'$ is
constant, where $\pi':X_0'\dashrightarrow C'$ is the cone
projection. Thus $\varphi(\C)\subseteq l$, where $l$ is a
generator of the cone $X_0'$.

By Proposition 1, $\varphi(\C)$ meets the double curve $D=X_0'\cap
X_0''$  of $X_0$ only in points of $D\cap X_{\infty}$. The
projection $\pi':D\to C'$ has degree $d''=4$ and simple
ramifications. Hence every fiber of
 $\pi'|D$ contains at least 3 points. A general
octic $X_{\infty}$ does not meet the ramification fibers of $\pi':
D\to C'$ and crosses $D$ passing through just one point of the
corresponding fiber of $\pi'|D$. Therefore $D \backslash
X_{\infty}$ contains at least 3 points of $l$. According to the
Little Picard Theorem, $\varphi: \C\to l\backslash (D \backslash
X_{\infty})$ is constant, a contradiction.

\vskip 0.2in

\noindent The Degeneration Principle of Proposition 1 can be
combined with the following one.

\vskip 0.2in

\noindent {\bf PROPOSITION 2 (Z$'$ 07$'$ \cite{Za})} {\it Let
$(X_t)_{t\in\PP^1}$ be a pencil of hypersurfaces in $\PP^{n+1}$
generated by two hypersurfaces $X_0$ and $X_\infty$ of the same
degree $d\ge 5$, where $X_0=kQ$ with $k\ge 2$ for some
hypersurface $Q\subseteq \PP^{n+1}$, and
$X_\infty=\bigcup_{i=1}^{d} H_{a_i}$, $a_1,\ldots,a_d\in\PP^1$, is
a union of $d$ distinct hyperplanes from a pencil
$(H_a)_{a\in\PP^1}$. If a sequence of Brody curves $\varphi_n:
\C\to X_{\varepsilon_n}$, where $\varepsilon_n\to 0$, converges to
a Brody curve $\varphi: \C\to X_0$, then $\varphi(\C)\subseteq
X_0\cap H_a$ for some $a\in\PP^1\,.$}

\vskip 0.2in

\noindent {\bf EXAMPLES} Given a pencil of planes $(H_a)$ in
$\PP^3$, using Proposition 2 one can deform  \\$\bullet$ $X_0=5Q$,
where $Q\subseteq\PP^3$ is a plane, \\$\bullet$ a triple quadric
$X_0=3Q\subseteq\PP^3$, or \\ $\bullet$ a double cubic, quartic,
etc.
$X_0=2Q\subseteq\PP^3$ \\
to an irreducible surface $X_\varepsilon\in\langle
X_0,X_\infty\rangle$ of the same degree $d$, where as before
$X_\infty=\bigcup_{i=1}^{d} H_{a_i}$, so that every limiting Brody
curve $\varphi:\C\to X_0$ is contained in a section $X_0\cap H_a$
for some $a\in\PP^1$.

\vskip 0.2in

\noindent The famous Bogomolov-Green-Griffiths-Lang Conjecture on
strong algebraic degeneracy (see e.g., \cite{BO,GG}) suggests that
every surface $S$ of general type possesses only finite number of
rational and elliptic curves and, moreover, the image of any
nonconstant entire curve $\varphi:\C\to S$ is contained in one of
them. In particular, this should hold for any smooth surface
$S\subseteq\PP^3$ of degree $\ge 5$, which fits the Kobayashi
Conjecture. Indeed, by Clemens-Xu-Voisin's Theorem, a general
smooth surface $S\subseteq\PP^3$ of degree $\ge 5$ does not
contain rational or elliptic curves, hence should be hyperbolic.
Anyhow, the  deformation method leads to the following result,
which is an immediate consequence of Proposition 2.

\vskip 0.2in

\noindent {\bf COROLLARY} {\it Let $S\subseteq\PP^3$ be a surface
and $Z\subset S$ be a curve such that the image of any nonconstant
entire curve $\varphi:\C\to S$ is contained in $Z$ \footnote{The
latter holds, for instance, if $S$ is hyperbolic modulo $Z$.}. Let
$X_\infty$ be the union of $d=2\deg S$ planes from a general
pencil of planes in $\PP^3$. Then any small enough linear
deformation $X_\varepsilon$ of $X_0=2S$ in direction of $X_\infty$
is hyperbolic.}

\vskip 0.2in

\noindent Along the same lines, Proposition 2 applies in the
following setting.

\vskip 0.2in

\noindent {\bf EXAMPLE} Let us take for $X_0$ a double cone in
$\PP^3$ over a plane hyperbolic curve of degree $d\ge 4$, and  for
$X_\infty$ a union of $2d$ distinct planes from a general pencil
$(H_a)$. Then small deformations $X_\varepsilon$ of $X_0$ in
direction of $X_\infty$ provide examples of hyperbolic surfaces of
any even degree $2d\ge 8$. For $d=4$ the latter surfaces can be
given by equation (\ref{1}) in suitable coordinates. Hence these
are actually the Duval-Fujimoto examples \cite{Du2, Fu}.

\vskip 0.2in

\noindent A nice construction due to J.\ Duval $'$04 \cite{Du} of
a hyperbolic sextic $X_{\varepsilon}\subseteq \PP^3$ uses the
deformation method iteratively in 5 steps, so that
$\varepsilon=(\varepsilon_1,\ldots,\varepsilon_5)$ has 5
subsequently small enough components. Hence $X_{\varepsilon}$ vary
within a 5-dimensional linear system; however the deformation of
$X_{0}$ to $X_{\varepsilon}$ neither is linear nor very generic.
It was suggested in \cite{SZ1} that the union of 6 general planes
in $\PP^3$ admits a general small linear deformation to an
irreducible hyperbolic sextic surface.

\vskip 0.2in

Let us finally turn to the Kobayashi problem on hyperbolicity of
complements of general hypersurfaces. By virtue of
Kiernan-Kobayashi-M.\ Green's version of Borel's Lemma, the
complement $\PP^n\setminus L$ of the union $L=\bigcup_{i=1}^{2n+1}
L_i$ of ${2n+1}$ hyperplanes in $\PP^n$ in general position is
Kobayashi hyperbolic. In particular, this applies to the union $l$
of 5 lines in $\PP^2$ in general position. Moreover \cite{Za1} $l$
can be deformed to a smooth quintic curve with hyperbolic
complement via a small deformation. This deformation proceeds in 5
steps and neither is linear nor very generic. So the following
question arises.

\medskip

 {\bf Question.}   Let $L$ ($M$) stands for the
union of $2n+1$ ($2n-1$, respectively) hyperplanes in $\PP^n$ in
general position. Is the complement of a general small linear
deformation of $L$ Kobayashi hyperbolic? Is a general small linear
deformation of $M$ Kobayashi hyperbolic? In particular, does the
union of 5 lines in $\PP^2$ (of 5 planes in  $\PP^3$) in general
position admit a general small linear deformation to an
irreducible quintic curve with hyperbolic complement (to a
hyperbolic quintic surface, respectively)?


\vskip 0.2in

\begin{thebibliography}{MmMm}

\bibitem{BO} Bogomolov F., De Oliveira B.\ Hyperbolicity of nodal
hypersurfaces. J.\ Reine Angew.\ Math.\ 596 (2006), 89--101.

\bibitem{CZ} Ciliberto C., Zaidenberg M.\
3-fold symmetric products of curves as hyperbolic hypersurfaces in
$\PP^4$. Intern.\ J.\ Math.\ 14 (2003), 413--436.

\bibitem{DEG} Demailly J.-P., El Goul J.\ Hyperbolicity of generic
surfaces of high degree in projective 3-space. Amer.\ J.\ Math.\
122 (2000), 515--546.

\bibitem{Du} Duval J.\ Une sextique hyperbolique dans
$\PP^3(\C)$. Math.\ Ann.\  330  (2004), 473--476.

\bibitem{Du2} Duval J.\ Letter to J.-P. Demailly, October 30, 1999
(unpublished).

\bibitem{Fu} Fujimoto H.\
A family of hyperbolic hypersurfaces in the complex projective
space. The Chuang special issue. Complex Variables Theory Appl.\
43 (2001), 273--283.

\bibitem{GG} Green M., Griffiths Ph.\
Two applications of algebraic geometry to entire holomorphic
mappings. The Chern Symposium 1979, 41--74, Springer, New
York-Berlin, 1980.

\bibitem{KZ} Kaliman S., Zaidenberg M.\
Non-hyperbolic complex spaces with hyperbolic normalization.
Proc.\ Amer.\ Math.\ Soc.\ 129 (2001), 1391--1393.


\bibitem{MQ} McQuillan M.\
Holomorphic curves on hyperplane sections of 3-folds. Geom.\
Funct.\ Anal.\ 9 (1999), 370--392.

\bibitem{Ro} Rousseau B.\ Equation diff\'erentielles sur
les hypersurfaces de $\PP^4$. J.\ Math\'em.\ Pure Appl.\ 86
(2006), 322--341.

\bibitem{SZ} Shiffman B., Zaidenberg M.\
New examples of Kobayashi hyperbolic surfaces in $\Bbb C\rm P\sp
3$. (Russian)  Funktsional. Anal. i Prilozhen.  39  (2005),
90--94; English translation in  Funct.\ Anal.\ Appl.\  39  (2005),
76--79.

\bibitem{SZ1} Shiffman B., Zaidenberg M.\
Constructing low degree hyperbolic surfaces in $\PP^3$. Special
issue for S. S. Chern. Houston J.\ Math.\ 28 (2002), 377--388.

\bibitem{SZ2} Shiffman B., Zaidenberg M.\
Hyperbolic hypersurfaces in $\PP^n$ of Fermat-Waring type. Proc.\
Amer.\ Math.\ Soc.\ 130 (2002), 2031--2035.

\bibitem{SZ3} Shiffman B., Zaidenberg M.\
Two classes of hyperbolic surfaces in $\PP^3$. International J.\
Math.\ 11 (2000), 65--101.

\bibitem{Si} Siu Y.-T.\ Hyperbolicity in Complex Geometry, in:
The legacy of Niels Henric Abel, Springer-Verlag, Berlin, 2004,
543--566.

\bibitem{Za} Zaidenberg M.\ Hyperbolicity of general deformations.
Preprint MPIM 106 (2007), 9p.

\bibitem{Za1} Zaidenberg M.\
Stability of hyperbolic embeddedness and construction of examples.
(Russian) Matem. Sbornik 135 (177) (1988), 361--372; English
translation in Math.\ USSR Sbornik 63 (1989), 351--361.

\end{thebibliography}
\end{document}